\documentstyle[12pt]{article}
\newtheorem{Lemma}{Lemma}
\newtheorem{Theorem}{Theorem}
\newtheorem{Prop}{Proposition}

\newtheorem{Cor}{Corollary}

\newcommand{\pf}{\medskip\noindent{\sc Proof: }}
\newcommand{\qed}{$\Box$}
\newcommand{\rTo}{\,\longrightarrow\,}

\newcommand{\DS}{\displaystyle}
\begin{document}
\title{The gamma function and a certain sequence of differences}
\author{David E. Radford \\
Department of Mathematics, Statistics \\
and Computer Science (m/c 249)    \\
851 South Morgan Street   \\
University of Illinois at Chicago\\
Chicago, Illinois 60607-7045
 }
\maketitle
\date{}
\begin{center}{\small Dedicated to the memory of Larry Lambe}\end{center}
\begin{abstract}
{\small  \rm We study a sequence of differences related to the problem of finding the smallest factorial $n!$ greater than or equal to $a^n$, where $a > 1$, using the gamma function.}
\end{abstract}
\let\thefootnote\relax\nonumber
\section{Introduction}\label{SectionGoal}
Let $a > 1$. Then $a^n \leq n!$ for some $n \geq 1$. Let $n_a$ be the smallest such positive integer. Then $n_a$ is the unique solution $n$ to $n \geq 2$ and $a$ belongs to the half open interval $(\sqrt[n-1]{(n-1)!}\,, \; \sqrt[n]{n!}\;]$. An extensive discussion of this problem is found in \cite{Radford}.

Let $n \geq 1$. In this paper we study the sequence of differences $\sigma_1, \sigma_2, \sigma_3, \ldots$, where $\sigma_n = \sqrt[n+1]{(n+1)!} - \sqrt[n]{n!}\;$, the sequence of lengths of the half open intervals, by way of the gamma function $\Gamma(x)$.

Let $G(x) = \Gamma(x+1)^{1/x}$, where $x > 0$. Then $\sigma_n = G(n+1) - G(n) = G'(x)$ for some $n < x < n+1$. Our study of the sequence is based on results we obtain about $G(x)$ and its first and second derivatives which are developed extensively in this paper. Our main result is Theorem \ref{TheoremMain} which states that the sequence $\sigma_1, \sigma_2, \sigma_3, \ldots \;$ is strictly decreasing and that $\DS{\lim_{n \rTo \infty}\sigma_n \frac{1}{e}}$. This theorem plays a central role in \cite{Radford}.

The main problem addressed in this paper is determining when $G'(x)$ is strictly decreasing. This we do by investigating when $G''(x)$ is negative. Section \ref{SecSequence} sets up the machinery for this investigation, particularly the other functions used in the process. Sections \ref{SecXFPrimeX}--\ref{SectionXGXDX2} are devoted to the technical details about these functions.

The main results are found in Section \ref{SecMainResults}. In addition to Theorem \ref{TheoremMain}, which is about the sequence $\sigma_1, \sigma_2, \sigma_3, \ldots \;$, Theorem \ref{TheoremSSequence} is about a closely related sequence $S_1, S_2, S_3, \ldots \;$, where $S_n = e\sigma_n$ for $n \geq 1$.

Let $A(x) = eG'(x)$, where $x > 0$, and let $n \geq 1$. Then $S_n = A(x)$ for some $n < x < n+1$.  Since $\lim_{x \rTo \infty}A(x) = 1$, it is natural to study $A(x)$ in terms of $a(x) = (A(x)-1)x$ as $\DS{A(x) = 1 + \frac{a(x)}{x}}$. The statements of Theorem \ref{TheoremMain} translate to statements about $S_1, S_2, S_3, \ldots$ in Theorem \ref{TheoremSSequence} which also states that $\DS{1 + \frac{a(n+1)}{n+1} < S_n < 1 + \frac{a(n)}{n}}$ for $n \geq 18$ and that $\DS{\lim_{n \rTo \infty} a(n) = \frac{1}{2}}$. We remark here that many of the functions $F(x)$ of interest in \cite{Radford} satisfy $\lim_{x \rTo \infty} F(x) = 1$ and are treated as $\DS{F(x) = 1 + \frac{f(x)}{x}}$.
%
%
\section{$\Gamma (x)$ and the sequence $\sigma_1, \sigma_2, \sigma_3, \ldots$}\label{SecSequence}
Let $\Gamma(x)$ denote the gamma function which is infinitely differentiable and positive valued on $(0, \infty)$. Let $x > 0$. Since $\Gamma (x+1) = x\Gamma(x)$ and $\Gamma(1) = 1$ it follows that $\Gamma(n) = (n-1)!$ for integers $n \geq 1$. See \S{1} and \S{2} of \cite{Artin} for the properties of the gamma function we use.

Let $G(x) = \Gamma(x+1)^{1/x} = e^{f(x)}$, where $\DS{f(x) = \frac{\ln \Gamma(x+1)}{x}}$. Then $G(x)$ and $f(x)$ are infinitely differentiable and $G(x)$ is positive valued on $(0, \infty)$. For $n \geq 1$, observe that
\begin{equation}\label{EqSigmaN}
\sigma_n = \sqrt[n+1]{(n+1)!} - \sqrt[n]{n!} = G(n+1) - G(n) = G'(x)
\end{equation}
for some $n < x < n+1$ by the Mean Value Theorem. We will show that $G'(x)$ is strictly decreasing eventually by investigating when $G''(x)$ is negative. Our results enable us to conclude that $\sigma_1, \sigma_2, \sigma_3, \ldots \;$ is strictly decreasing sequence.

We first observe that
\begin{equation}\label{EqNegative}
G''(x) < 0 \;\; \mbox{if and only if} \;\; f''(x) + f'(x)^2 < 0
\end{equation}
since
\begin{equation}\label{EqBasicGPrimePrime}
G''(x) = (f''(x) + f'(x)^2)G(x).
\end{equation}
As $\Gamma(x+1) = x\Gamma(x)$, we have
\begin{equation}\label{EqXFX}
xf(x) = \ln x + \ln \Gamma(x).
\end{equation}
Now
$$
\DS{\ln \Gamma(x) = -Cx - \ln x + \sum_{n = 1}^\infty \left(\frac{x}{n} - \ln(1 + \frac{x}{n})\right)},
$$
where $C$ is the Euler-Mascheroni constant, by (2.9) of \cite{Artin}. Therefore
$$f(x) = -C + \DS{\sum_{n=1}^\infty \left(\frac{1}{n} - \frac{1}{x}\ln(\DS{1 + \frac{x}{n}})\right)}.$$
As
$$
(\ln \Gamma)'(x) = -C - \frac{1}{x} + \sum_{n = 1}^\infty \left(\frac{1}{n} - \frac{1}{x+n}\right)
$$
by (2.10) of \cite{Artin}, and
\begin{equation}\label{EqXFXPrimeEq}
\DS{xf'(x) = \frac{1}{x} + (\ln \Gamma)'(x) - f(x)}
\end{equation}
by (\ref{EqXFX}), the equation
\begin{equation}\label{EqXFXPrime}
xf'(x) = \sum_{n = 1}^\infty \left(\frac{1}{x} \ln(1 + \frac{x}{n}) - \frac{1}{x+n}\right)
\end{equation}
holds. We can look at the terms of the series in the expressions for $\ln \Gamma(x)$ and $xf'(x)$ in light of estimates $\DS{\frac{a}{1+a} < \ln (1 + a) < a}$ for $\ln (1+a)$, where $a > 0$. The typical summand in (\ref{EqXFXPrime}) is $\DS{\frac{1}{x}\left(\ln(1+a) - \frac{a}{1+a}\right)}$, where $\DS{a = \frac{x}{n}}$. As a consequence $xf'(x) > 0$. Now $\DS{xf''(x) = -\frac{1}{x^2} + (\ln \Gamma)''(x) - 2f'(x)}$ by (\ref{EqXFXPrimeEq}). Since $\DS{(\ln \Gamma)''(x) = \sum_{n = 0}^\infty \frac{1}{(x+n)^2}}$ by (2.11) of \cite{Artin}, we have
\begin{equation}\label{EqX2FPrimePrime}
x^2(f''(x) + f'(x)^2) = \sum_{n = 1}^\infty \frac{x}{(x+n)^2} + (1- xf'(x))^2 - 1.
\end{equation}
This equation and (\ref{EqBasicGPrimePrime}) imply
\begin{equation}\label{EqMain}
G''(x) < 0 \;\; \mbox{if and only if} \;\; \DS{\sum_{n = 1}^\infty \frac{x}{(x+n)^2} + (1 - xf'(x))^2 < 1}.
\end{equation}
The equation of (\ref{EqXFXPrime}) and the two subsequent ones imply
\begin{equation}\label{EqX2F2Primes}
x^2f''(x) = 2\left(\DS{\sum_{n = 1}^\infty \left(\frac{3x+2n}{2(x+n)^2} - \frac{1}{x} \ln (1 + \frac{x}{n})\right)}\right).
\end{equation}
Fix $x > 0$ and set $D(t) = \DS{\frac{1}{x}}\ln (1 + \frac{x}{t}) - \frac{3x + 2t}{2(x+t)^2}$ for $t > 0$. Then $\lim_{t \rTo \infty}D(t) = 0$ and $D'(t) < 0$. Therefore $D(t)$ is a positive valued strictly decreasing function. Note $x^2f''(x) = -2(\sum_{n = 1}^\infty D(n))$.

An antiderivative of $D(t)$ is $\DS{\frac{t}{x}\ln (1+ \frac{x}{t}) - \frac{t}{2(x+t)}}$. As a consequence $\DS{\int_1^\infty D(t)dt} = \DS{\frac{1}{2} - \left(\frac{1}{x}\ln (1+x) - \frac{1}{2(x+1)} \right)}$.
Using the integral test one can easily see that $\DS{\lim_{x \rTo \infty}\left(\sum_{n = 1}^\infty \left(\frac{3x+2n}{2(x+n)^2} - \frac{1}{x} \ln (1 + \frac{x}{n})\right)\right) = -\frac{1}{2}}$. We have shown $x^2f''(x) < 0$ and
\begin{equation}\label{EqX2F2Primes2}
\lim_{x \rTo \infty} x^2f''(x) = -1.
\end{equation}
Since $-f''(x) > 0$ for all $x > 0$, it follows by the corollary to  Theorem 1.4 of \cite{Artin} that $-f(x)$ is a convex function on $(0, \infty)$.

Note that $D(1) > 0$ means $\DS{\frac{a}{1+a} < \frac{3a^2 + 2a}{2(1+a)^2} < \ln (1+a)}$ for all $a > 0$, which gives a better underestimate for $\ln (1+a)$ than its usual underestimate $\DS{\frac{a}{a+1} < \ln (1+a)}$. Using (\ref{EqXFXPrime}) and (\ref{EqX2F2Primes}), it is easy to see that $\DS{f'(x) > -\frac{x}{2}f''(x)}$ for all $x > 0$. One way to prove this inequality is to think of the terms of the series involved in terms of the usual and better under estimates for $\ln(1+a)$.

We will consider the two left hand terms in the second inequality of (\ref{EqMain}) separately. As for the first, we set $g(x) = \DS{\sum_{n = 1}^\infty \frac{x}{(x+n)^2}}$. Observe that $g(x) = x^2f''(x) + 2xf'(x)$ by (\ref{EqX2FPrimePrime}). In light of the better underestimate for $\ln (1+a)$, it is interesting to compute $g(x)$ using (\ref{EqXFXPrime}) and (\ref{EqX2F2Primes}).

Fix $x > 0$ and set $\DS{D(t) = \frac{x}{(x+t)^2}}$ where $t > 0$. Observe that $D(t)$ is a positive valued strictly decreasing function and $\lim_{t \rTo \infty} D(t) = 0$. As $\DS{\int_1^\infty D(t)dt = \frac{x}{x+1}}$, and $\DS{\frac{x}{x+1} + D(1) = 1 - \frac{1}{(x+1)^2}}$, by the integral test for series
\begin{equation}\label{EqFirstSeries}
\frac{x}{x+1} < g(x) < 1 - \frac{1}{(x+1)^2}.
\end{equation}
As a result $\lim_{x \rTo \infty}g(x) = 1$. When $m$ is positive integer, $g(x)$ has a particularly nice form. Since $\zeta(2) = \DS{\sum_{n = 1}^\infty \frac{1}{n^2}} = \frac{\pi^2}{6}$, where $\zeta(x)$ is the Riemann zeta function,
\begin{equation}\label{EqGM}
g(m) = m(\frac{\pi^2}{6} - \sum_{n = 1}^m \frac{1}{n^2}).
\end{equation}
We consider the function $xf'(x)$ in the next section.
\section{$xf'(x)$ and applications to $G(x)$}\label{SecXFPrimeX}
Throughout $x > 0$. Let
$$
h(x) = xf'(x) = \sum_{n = 1}^\infty \left(\frac{1}{x} \ln(1 + \frac{x}{n}) - \frac{1}{x+n}\right)
$$
and $d(x) = 1 - h(x)$. We begin our discussion by considering the infinite series $h(x)$.

Fix $x$ and set $D(t) = \DS{\frac{1}{x}\ln(1 + \frac{x}{t}) - \frac{1}{x+t}}$ for $t > 0$. Using the usual underestimate for $\ln (1+a)$ one sees that $\DS{\frac{x}{x+t} < \ln (1 + \frac{x}{t})}$. Thus $D(t) > 0$. It is easy to see that $\lim_{t \rTo \infty}D(t) = 0$. Since $D'(t) < 0$, we conclude $D(t)$ is a strictly decreasing function. An antiderivative of $D(t)$ is $\DS{\left(\frac{t}{x} \right)\ln \left(1 + \frac{x}{t}\right)}$.  Thus we compute $\DS{\int_1^\infty D(t)dt = 1 - \frac{\ln (1 + x)}{x}}$. Using the integral test for series we deduce
\begin{equation}\label{EqHXEstimates}
1 - \frac{\ln (1+x)}{x} < h(x) <  \frac{x}{x+1}.
\end{equation}
Therefore
\begin{equation}\label{EqLimitXFPtrimeX}
\lim_{x \rTo \infty} h(x) = 1 \;\; \mbox{and} \;\; \lim_{x \rTo \infty} d(x) = 0,
\end{equation}
\begin{equation}\label{EqSecondSeriesEstimates1}
\frac{1}{x+1} < d(x) < \frac{\ln (1+x)}{x},
\end{equation}
and thus
\begin{equation}\label{EqSecondSeriesEstimates2}
\frac{1}{(x+1)^2} < d(x)^2 < \left(\frac{\ln (1+x)}{x}\right)^2.
\end{equation}
The first limit of (\ref{EqLimitXFPtrimeX}) follows by (\ref{EqX2F2Primes2}) and the discussion of $g(x)$ in Section \ref{SecSequence} as well.

Our integral test approximations give us $L(x) < g(x) + d(x)^2 < R(x)$, where $\DS{L(x) = 1 - \frac{x}{(x+1)^2}}$ and where $\DS{R(x) = 1 + \left(\frac{\ln (1+x)}{x}\right)^2 - \frac{1}{(x+1)^2}}$. See (\ref{EqFirstSeries}) and (\ref{EqSecondSeriesEstimates2}). Since $\DS{\frac{1}{x+1} < \frac{\ln(1 + x)}{x}}$, these approximations are not sufficient to give the second inequality of (\ref{EqMain}).

We return to (\ref{EqLimitXFPtrimeX}) for a moment. In light of (\ref{EqX2F2Primes2}) we have
\begin{equation}\label{EqMain2}
\lim_{x \rTo \infty} x^2(f''(x) + f'(x)^2) = 0 = \lim_{x \rTo \infty} (f''(x) + f'(x)^2).
\end{equation}

We now compute $h(x)$ when $x = m$ is a positive integer. In this case $h(x)$ is expressed in terms of the Euler-Mascheroni constant and terms of the natural sequence which defines it. Our calculations will lead to other properties of $G(x)$ and its first and second derivatives.

Let $N \geq 1$ and set
$$
h_N(x) = \sum_{n = 1}^N \left(\frac{1}{x} \ln(1 + \frac{x}{n}) - \frac{1}{x+n}\right)  = \frac{1}{x}\left(\sum_{n = 1}^N\left(\ln(1 + \frac{x}{n}) - \frac{x}{x+n}\right)\right).
$$
Let $m \geq 1$. Then
$$
h_N(m) = \frac{1}{m}\left(\sum_{n = 1}^N\left(\ln (m+n) - \ln n\right) - m\left(\sum_{n = 1}^N\frac{1}{m+n}\right)\right).
$$
To proceed we use a very simple summation formula.
\begin{Lemma}\label{LemmaSums}
Let $a_1, a_2, a_3, \ldots$ be any sequence of real numbers and suppose $N \geq m \geq 1$. Then  $\DS{\sum_{n=1}^N (a_{m+n} - a_n) = \sum_{\ell = 1}^m (a_{N+ \ell} - a_\ell)}$.
\end{Lemma}

\pf
The lemma follows from $\sum_{n=1}^N a_{m+n} = \sum_{\ell=1}^{N-m}a_{m+\ell} + \sum_{\ell = 1}^ma_{N+\ell}$ and from $\sum_{n =1}^N a_n = \sum_{\ell = 1}^m a_\ell + \sum_{\ell = 1}^{N-m}a_{m+\ell}$.
\qed
\medskip

Assume further that $m \leq N$. In light of the preceding lemma
\begin{equation}\label{EqHNXRewrite}
h_N(m) = \frac{1}{m}\left(\sum_{\ell = 1}^m\left(\ln (N+ \ell) - \ln \ell - \sum_{n = 1}^N \frac{1}{m+n}\right)  \right).
\end{equation}

We consider each summand separately.
$$
\ln (N+\ell) - \ln \ell - \sum_{n = 1}^N\frac{1}{m+n} = \left(\ln (N + \ell) - \sum_{j = 1}^{N+m}\frac{1}{j}\right) - \left(\ln \ell- \sum_{j = 1}^m\frac{1}{j}\right)
$$
Let $C_n = 1 + 1/2 + 1/3 + \cdots  + 1/n - \ln n$ for $n \geq 1$. Then $C_1, C_2, C_3, \ldots \;$ is a strictly descending sequence of positive numbers which converges to $C$. See page 16 of \cite{Artin}. Therefore
$$
\ln (N+\ell) - \ln \ell - \sum_{n = 1}^N\frac{1}{m+n} = -C_{N+\ell} + C_\ell- \sum_{j = \ell +1}^m\frac{1}{N+j} + \sum_{j = \ell+1}^m \frac{1}{j}.
$$
Since $h(m) = \lim_{N \rTo \infty}h_N(m)$, by (\ref{EqHNXRewrite}) we now have
$$
h(m) = \frac{1}{m}\left(-mC + \sum_{\ell = 1}^mC_\ell + \sum_{\ell = 1}^m\left(\sum_{j = \ell+1}^m \frac{1}{j} \right)\right).
$$
Now $\DS{\sum_{\ell = 1}^m\left(\sum_{j = \ell+1}^m \frac{1}{j} \right)}$ $ = \DS{m - \sum_{j = 1}^m \frac{1}{j}}$ $= m - C_m - \ln m$. The first equation easily follows by induction on $m$. Therefore
\begin{equation}\label{EqHNXFinal}
h(m) = 1 - C  + \frac{C_1 + \cdots + C_{m-1} - \ln m}{m}
\end{equation}
for $m \geq 1$, where $C_1 + \cdots + C_{m-1} = 0$ when $m = 1$ by convention. As a result
\begin{equation}\label{EqOneMinusHM}
\frac{1}{m+1} < d(m) = C  - \frac{C_1 + \cdots + C_{m-1} - \ln m}{m} < \frac{\ln (1 + m)}{m}
\end{equation}
by (\ref{EqSecondSeriesEstimates1}). We will simplify the expressions for $d(m)$ and $h(m)$.
\begin{Lemma}\label{LemmaSimplifySums}
$\DS{\sum_{\ell = 1}^mC_\ell =  (m+1)(C_{m+1} + \ln (m + 1) - 1) - \ln (m!)}$ for $m \geq 0$.
\end{Lemma}
\pf
$\DS{\sum_{\ell = 1}^m\left(\sum_{j = 1}^\ell\frac{1}{j}\right) + \sum_{\ell = 1}^m\left(\sum_{j = \ell+1}^m\frac{1}{j}\right) = \sum_{\ell = 1}^m\left(\sum_{j = 1}^m\frac{1}{j}\right) = m\left(\sum_{j = 1}^m\frac{1}{j}  \right)   }$ and

$\DS{\sum_{\ell = 1}^m\left(\sum_{j = \ell+1}^m \frac{1}{j} \right) = m - \sum_{j = 1}^m \frac{1}{j}}$ imply $\DS{\sum_{\ell = 1}^m\left(\sum_{j = 1}^\ell\frac{1}{j}\right) = (m+1)\left(\sum_{j = 1}^m\frac{1}{j}\right) - m}$
\medskip

\noindent
from which the lemma easily follows. \qed
\medskip

For $m \geq 1$ let $\DS{D_m = 1 + \frac{1}{2} + \cdots + \frac{1}{m} - \ln (\sqrt[m]{m!})}$. By the preceding lemma and (\ref{EqOneMinusHM}) we conclude:
\begin{Cor}\label{CorHM}
$\DS{d(m) = C + 1 - D_m}$ and $\DS{h(m) =  D_m - C}$ for all $m \geq 1$. \qed
\end{Cor}
\medskip

Note that (\ref{EqLimitXFPtrimeX}) and Corollary \ref{CorHM} imply that $\lim_{m \rTo \infty} D_m = C + 1$. We will examine the sequence $D_1, D_2, D_3, \ldots \;$ in more detail. Observe that $\DS{D_m = C_m + \ln m - \ln (\sqrt[m]{m!})}$ $\DS{ = C_m + \ln (\frac{m}{\sqrt[m]{m!}})}$ and hence $\DS{\lim_{m \rTo \infty} \frac{m}{\sqrt[m]{m!}} = e}$, or equivalently  $\DS{\lim_{m \rTo \infty} \frac{\sqrt[m]{m!}}{m} = \frac{1}{e}}$, follows.

Now $\DS{D_{m+1} = D_m + \frac{1}{m+1}\left(1 + \ln \left(\frac{\sqrt[m]{m!}}{m+1}\right)\right)}$. Recall from Section \ref{SecSequence} that $\sqrt[m]{m!} = G(m) = e^{f(m)}$. Thus $\DS{\ln \left(\frac{\sqrt[m]{m!}}{m+1}\right)} = \DS{\ln \left(\frac{G(m)}{m+1}\right)} = f(m) - \ln (m+1)$.

Let $D(x) = f(x) - \ln (x+1) = \DS{\ln \left(\frac{G(x)}{x+1}\right)}$. The second inequality of (\ref{EqHXEstimates}) can be expressed $xD'(x) < 0$. This means that $D(x)$, and hence $\DS{\frac{G(x)}{x+1}}$, is a strictly decreasing function. By virtue of our last limit calculation  $\DS{\lim_{m \rTo \infty} \ln \left(\frac{\sqrt[m]{m!}}{m+1}\right) = -1}$. We have shown that $D_1, D_2, D_3, \ldots \:$ is a strictly increasing sequence.

Since $D(x)$ is a strictly decreasing function and $\lim_{m \rTo \infty}D(m) = -1$, it follows that $\lim_{x \rTo \infty}D(x) = -1$. We have observed that $xf'(x) > 0$ in the discussion following by (\ref{EqXFXPrime}). Thus $\DS{(x+1)f'(x) =  \left(\frac{x+1}{x}\right)xf'(x) > 0}$ and $\lim_{x \rTo \infty} (x+1)f'(x) = 1$ by (\ref{EqLimitXFPtrimeX}). As $G'(x) = f'(x)G(x) = e^{E(x)}$, where $E(x) = D(x) + \ln ((x+1)f'(x))$, and $\lim_{x \rTo \infty}E(x) = -1 = \lim_{x \rTo \infty}D(x)$, we conclude that
\begin{equation}\label{EqGPrimeLimit}
\lim_{x \rTo \infty}G'(x) = \frac{1}{e} = \lim_{x \rTo \infty} \frac{G(x)}{x}.
\end{equation}
Recall that $G''(x) = (f''(x) + f'(x)^2)G(x)$ by (\ref{EqBasicGPrimePrime}). Thus
\begin{equation}\label{EqGPrimePrimeLimit}
\lim_{x \rTo \infty}xG''(x) = 0
\end{equation}
follows by (\ref{EqMain2}) and (\ref{EqGPrimeLimit}). We end this section with
\begin{equation}\label{EqGDerivativeLeftRight}
\frac{G(x)}{x+1}R(x) \leq G'(x) \leq \frac{G(x)}{x+1},
\end{equation}
where $\DS{R(x) = \left(\frac{1+x}{x^2}\right)(x - \ln (1 + x))}$. It is easy to see that $0 < R(x) < 1$, that $R(x)$ is a strictly increasing function, and that $\lim_{x \rTo \infty}R(x) = 1$.
To establish the second inequality of (\ref{EqGDerivativeLeftRight}), we let $h > 0$. Since $\DS{\frac{G(x)}{x+1}}$ is a strictly decreasing function, we have that $\DS{\frac{G(x+h)}{x+h+1} < \frac{G(x)}{x+1}}$ from which $\DS{\frac{G(x+h) - G(x)}{h} < \frac{G(x)}{x+1}}$ follows. The second inequality is established.

Now let $F(x) = \DS{\frac{1}{x}\ln (1 + x) + \ln (1 + x)}$. Then the first inequality of  (\ref{EqHXEstimates}) can be written $xF'(x) < xf'(x)$. Thus $\DS{\ln \left(\frac{G(x)}{(x+1)\ell(x)}\right)}$ $= f(x) - F(x)$, hence $\DS{\frac{G(x)}{(x+1)\ell(x)}}$, is a strictly increasing function, where $\ell(x) = (1 + x)^{1/x}$. At this point we leave it to the reader to establish the first inequality of (\ref{EqGDerivativeLeftRight}) and assertions about $R(x)$, starting by mimicking our proof of the second inequality.
\section{$a(x)$ and $G''(x)$}\label{SectionAX}
Let $x > 0$ and set $A(x) = eG'(x)$. Now $\lim_{x \rTo \infty} A(x) = 1$ by (\ref{EqGPrimeLimit}). Let $a(x) = (A(x)-1)x$. Then $\DS{A(x) = 1 + \frac{a(x)}{x}}$. We find $a(x)$ is useful to associate with $A(x)$. As $a'(x) = A'(x)x + A(x) - 1 = e(xG''(x) + G'(x)) - 1$ we have
\begin{equation}\label{EqAPrimeLimit}
\lim_{x \rTo \infty} a'(x) = 0
\end{equation}
by (\ref{EqGPrimeLimit}) and (\ref{EqGPrimePrimeLimit}). To proceed we need estimates for a certain expression involving the derivative of the polygamma function $\psi (x) = (\ln \Gamma)'(x)$.
\begin{Lemma}\label{LemmaPsi}
Let $x> 0$. Then:
\begin{enumerate}
\item[{\rm (a)}] $\DS{\frac{1}{2} + \frac{x^2}{6(x + 1/14)^3} < x^2\psi'(x) - x < \frac{1}{2} + \frac{1}{6x}}$.
\item[{\rm (b)}] $\DS{\lim_{x \rTo \infty} (x^2\psi'(x) - x) = \frac{1}{2}}$.
\end{enumerate}
\end{Lemma}

\pf
Part (b) follows from part (a) which in turn follows from the two inequalities
$$
\DS{\frac{1}{x} + \frac{1}{2x^2} + \frac{1}{6(x + 1/14)^3} < \psi'(x) < \frac{1}{x} + \frac{1}{2x^2} + \frac{1}{6x^3}}
$$
described in Theorem 4 of \cite{Gordon}. \qed
\medskip

Since $\DS{g(x) = x\left(\psi'(x) - \frac{1}{x^2}\right)}$, it follows that $xg(x) - x = x^2\psi'(x) - x - 1$. We will need an equivalent formulation of the preceding lemma:
\begin{Cor}\label{CorXGXMinusX}
Let $x > 0$. Then:
\begin{enumerate}
\item[{\rm (a)}] $\DS{-\frac{1}{2} + \frac{x^2}{6(x + 1/14)^3} < xg(x) - x < -\frac{1}{2} + \frac{1}{6x}}$.
\item[{\rm (b)}] $\DS{\lim_{x \rTo \infty} (xg(x) - x) = -\frac{1}{2}}$.
\end{enumerate} \qed
\end{Cor}
\medskip

Now
\begin{eqnarray*}
x^2G''(x) & = & (g(x) + d(x)^2 - 1)G(x) \\
& = & (xg(x) - x + xd(x)^2)\frac{G(x)}{x}
\end{eqnarray*}
follows by (\ref{EqBasicGPrimePrime}) and (\ref{EqX2FPrimePrime}). Observe that $\DS{\lim_{x \rTo \infty}xd(x)^2 = 0}$ by (\ref{EqSecondSeriesEstimates2}) and recall that $\DS{\lim_{x \rTo \infty}\frac{G(x)}{x} = \frac{1}{{e}}}$ by (\ref{EqGPrimeLimit}). Hence
\begin{equation}\label{EqAXGPrimePrimeLimit}
\lim_{x \rTo \infty} x^2G''(x) = -\frac{1}{2e}.
\end{equation}
now follows by part (b) of Corollary \ref{CorXGXMinusX}. Since $\lim_{x \rTo \infty}A(x) = 1$, in light of (\ref{EqAXGPrimePrimeLimit}) we have
\begin{equation}\label{EqALimit}
\lim_{x \rTo \infty}a(x) = \frac{1}{2}
\end{equation}
by L'H\^{o}spital's Rule. Observe (\ref{EqAXGPrimePrimeLimit}) implies $G'(x)$ is eventually decreasing.
\section{$xg(x) - x + xd(x)^2$ revisited}\label{SectionXGXDX2}
By the calculation following the proof of Corollary \ref{CorXGXMinusX} we see that $G''(x) < 0$ if and only if $xg(x) - x +xd(x)^2 < 0$. Let $x > 0$. Then
$$
\DS{xg(x) - x +xd(x)^2 < -\frac{1}{2} + \frac{1}{6x} + \frac{(\ln(1+x))^2}{x}}
$$
by part (a) of Corollary \ref{CorXGXMinusX} and (\ref{EqSecondSeriesEstimates2}).  The right hand side of this inequality is negative from some point on. When $\DS{D(x) = \frac{3x-1}{6} - (\ln (1+x))^2}$ becomes and remains positive this is the case. Note that $\DS{D'(x) = \frac{1}{2} - 2 \left(\frac{\ln (1+x)}{1+x}\right)}$. Since $\DS{D''(x) = -2\left(\frac{1 - \ln (1+x)}{(1+x)^2}\right)}$, it follows $D'(x)$ is strictly decreasing on $(0, e-1]$ and is strictly increasing on $[e-1, \infty)$. Since $D'(e-1) < 0$ and $\lim_{x \rTo \infty}D'(x) = 1/2$, it follows that $D'(x)$ has a unique zero $a > e-1$, which is the unique solution to this inequality and $\DS{\frac{\ln(1+a)}{1+a} = \frac{1}{4}}$, and $D(x)$ is strictly increasing on $[a, \infty)$. As $\DS{\frac{\ln 9}{9} < \frac{1}{4} < \frac{\ln 8}{8}}$ it follows that $7 < a < 8$. As $D(17) < 0 < D(18)$, we have shown part (a) of the following:
\begin{Prop}\label{PropMainPrize}
For $G(x) = \Gamma(x+1)^{1/x}$ the following hold:
\begin{enumerate}
\item[{\rm (a)}] There exists $17 < c < 18$ such that $G''(x) < 0$ for $x > c$. Thus $G'(x)$ is strictly decreasing on $[c, \infty)$.
\item[{\rm (b)}] $\DS{\lim_{x \rTo \infty}xG''(x) = 0}$ and $\DS{\lim_{x \rTo \infty} x^2G''(x) = -\frac{1}{2e}}$.
\item[{\rm (c)}] $\DS{\lim_{x \rTo \infty}G'(x) = \frac{1}{e} = \lim_{x \rTo \infty} \frac{G(x)}{x}}$.
\end{enumerate} \qed
\end{Prop}
\medskip

Part (b) is a restatement of (\ref{EqGPrimePrimeLimit}) and (\ref{EqAXGPrimePrimeLimit}) respectively. Part (c) is a restatement of (\ref{EqGPrimeLimit}).

We note that $7.61316 < a < 7.61317$. Our $c$ of part (a), the unique solution to $c \geq a$ and $D(c) = 0$, satisfies $17.11650 < c < 17.11651$ as $D(17.11650) < 0 < D(17.11651)$.
\begin{Cor}\label{CorEuler}
$\DS{m\left(\frac{\pi^2}{6} - \sum_{n=1}^m\frac{1}{n^2}\right) + \left(C - \frac{1}{m}\left(\sum_{n = 1}^{m-1}C_n - \ln m\right)\right)^2} < 1$ for all $m \geq 1$.
\end{Cor}

\pf
First note that $G''(x) < 0$ if and only if $g(x) + d(x)^2) < 1$ by (\ref{EqBasicGPrimePrime}) and (\ref{EqX2FPrimePrime}) since $G(x) > 0$. When $m \geq 18$, the corollary follows (\ref{EqGM}), by (\ref{EqOneMinusHM}), and by part (a) of Proposition \ref{PropMainPrize}. The remaining cases follow by numerical calculation. \qed
\medskip

Note that $\DS{\left(C - \frac{1}{m}\left(\sum_{n = 1}^{m-1}C_n - \ln m\right)\right)^2}$ $= (D-D_m)^2$, where $D = C+1$, by Corollary \ref{CorHM}. We can easily restate Corollary \ref{CorEuler} in terms of the polygamma function $\psi'(x)$. That $\DS{g(x) = x\psi'(x) - \frac{1}{x}}$ follows by the calculation preceding Corollary \ref{CorXGXMinusX}. Thus:
\begin{Cor}\label{CorEulerPolyGamma}
$\DS{m\psi'(m) + \left(C - \frac{1}{m}\left(\sum_{n = 1}^{m-1}C_n - \ln m\right)\right)^2} < 1 + \frac{1}{m}$ holds for all $m \geq 1$. \qed
\end{Cor}
\section{The main results}\label{SecMainResults}
\begin{Theorem}\label{TheoremMain}
Let $\sigma_n = \sqrt[n+1]{(n+1)!} - \sqrt[n]{n!}$ for $n \geq 1$. Then:
\begin{enumerate}
\item[{\rm (a)}] $\lim_{n \rTo \infty}\sigma_n = \DS{\frac{1}{e}}$.
\item[{\rm (b)}] $0.42 > \sigma_1 > \sigma_2 > \sigma_3 > \cdots > \DS{\frac{1}{e}}$.
\end{enumerate}
\end{Theorem}

\pf
Let $n \geq 1$. Then $\sigma_n = G'(x)$ for some $n < x < n+1$ by (\ref{EqSigmaN}). Part (a) thus follows by (\ref{EqGPrimeLimit}). Since $G'(x)$ is strictly decreasing on $[18, \infty)$ by part (a) of Proposition \ref{PropMainPrize}, we have shown that $\sigma_{18} > \sigma_{19} > \sigma_{20} > \cdots \;$. Direct computation shows that $0.42 > \sigma_1 > \sigma_2 > \sigma_3 > \cdots > \sigma_{18}$. Thus the sequence $\sigma_1, \sigma_2, \sigma_3, \ldots \;$ is strictly decreasing and part (b) now follows by part (a).
\qed
\medskip

The sequence $S_1, S_2, S_3, \ldots \;$ is fundamental to \cite{Radford}, where $S_n = e\sigma_n$ for all $n \geq 1$. To study this sequence we use the functions $A(x)$ and $a(x)$ of Section \ref{SectionAX}. Let $x > 0$. Recall that $A(x) = eG'(x)$ and $\DS{A(x) = 1 + \frac{a(x)}{x}}$. Now $A(x)$ is strictly decreasing on $[18, \infty)$  by part (a) of Proposition \ref{PropMainPrize}. Let $n \geq 18$. As $S_n = A(x)$ for some $n < x < n+1$ by (\ref{EqSigmaN}); thus $A(n) > S_n > A(n+1)$ follows. These inequalities,Theorem \ref{TheoremMain}, and (\ref{EqALimit}) imply the following:
\begin{Theorem}\label{TheoremSSequence}
Let $S_n = e\sqrt[n+1]{(n+1)!} - e\sqrt[n]{n!} \;$ for $n \geq 1$. Then:
\begin{enumerate}
\item[{\rm (a)}] $\lim_{n \rTo \infty}S_n = 1$.
\item[{\rm (b)}] $1.15 > S_1 > S_2 > S_3 > \cdots > 1$.
\item[{\rm (c)}] $1 + \DS{\frac{a(n+1)}{n+1}} < S_n < 1 + \DS{\frac{a(n)}{n}}$ for all $n \geq 18$.
\item[{\rm (d)}] $\lim_{n \rTo \infty}a(n) = \DS{\frac{1}{2}}$.
\end{enumerate} \qed
\end{Theorem}
\medskip

In \cite{Radford} we consider over estimates for $S_n$ in detail. Our work based on the improvement on Sterling's approximation of factorials due to Robbins \cite{Robbins}. A few final comments based on parts (c) and (d) of Theorem \ref{TheoremSSequence}.

Part (c) can be expressed in a different manner as $\DS{S_n < 1 + \frac{a(n)}{n}}$ for $n = 18$ and $(S_n-1)n < a(n) < (S_{n-1}-1)n$ for $n > 18$. The last three inequalities imply $\DS{\lim_{n \rTo \infty}(S_n - 1)n = \frac{1}{2}}$. By part (b) there exists a sequence $a_1, a_2, a_3, \ldots \;$ such that $S_1 < 1 + a_1$ and $(S_n-1)n < a_n < (S_{n-1}-1)n$ for $n > 1$. For any such sequence $\DS{1 + \frac{a_{n+1}}{n+1} < S_n < 1 + \frac{a_n}{n}}$ for all $n \geq 1$ and $\DS{\lim_{n \rTo \infty}a_n = \frac{1}{2}}$.

\end{document}